\documentclass{amsart}

\usepackage{amsthm}
\usepackage[ruled,vlined,linesnumbered]{algorithm2e}
\usepackage[T1]{fontenc}
\usepackage{booktabs}
\usepackage{graphicx}
\usepackage{dcolumn}
\usepackage{caption}
\usepackage{url}
\usepackage{enumitem}

\newtheorem{theorem}{Theorem}[section]
\newtheorem{corollary}[theorem]{Corollary}

\newtheorem{lemma}[theorem]{Lemma}
\newtheorem{proposition}[theorem]{Proposition}

\theoremstyle{definition}
\newtheorem{definition}[theorem]{Definition}

\newcommand{\sign}{\operatorname{sign}} 

\title[Dynamic invisibility of implicit state-dependent delays]
{Dynamic invisibility of implicit and nested state-dependent delays in a class of externally forced delay equations}

\author{S. Bolduc-St-Aubin}
\address{Department of Mathematics, University of Auckland, Private Bag 92019, Auckland 1142, New Zealand}
\email{sbol082@aucklanduni.ac.nz}

\author{P. Subramanian}
\address{Department of Mathematics, University of Auckland, Private Bag 92019, Auckland 1142, New Zealand}

\author{B. Krauskopf}
\address{Department of Mathematics, University of Auckland, Private Bag 92019, Auckland 1142, New Zealand}

\date{October 2026}

\begin{document}

\begin{abstract}
Delayed feedback is a fundamental mechanism in many systems and is often modeled by delay differential equations (DDEs). While delays are frequently assumed to be constant, in realistic feedback processes they may depend on the state and even be defined implicitly or in a nested fashion by the feedback relation itself. 

In this work, we address the question of how such implicit and nested state-dependent delays affect the observed behavior of a class of externally forced systems with negative switched feedback. Our motivating example and principal case study is the piecewise-smooth limiting form of a conceptual DDE model for the El Ni{\~n}o--Southern Oscillation (ENSO) with a negative feedback term subject to implicit state-dependent delay. We prove, for a wider prototypical class of DDEs with switched negative feedback and more general forcing, that (subject to an explicit condition on the parameters) any finite level of nested state-dependent delay and the limiting implicit delay are dynamically invisible: they all generate exactly the same trajectory, namely, that of the corresponding constant-delay DDE. We also discuss how the delayed argument functions organize the temporal order of reactions (TOR) in this regime, and show that a more stringent condition is required to ensure TOR for nested state-dependent delays.
\end{abstract}

\maketitle

\section{Introduction}

Negative feedback, where the output counteracts changes, is a stabilization mechanism that plays a vital role in natural and man-made systems. Due to inevitable transmission and reaction/processing times, any such feedback is typically acting only after some time delay $\tau$. This leads to a description of such systems by delay differential equations (DDEs), in which the time derivative also depends on past values
of the state. Delays are often assumed to be constant, that is, remain unchanged as the system evolves. This is a helpful assumption because the respective model is then from the class of constant-delay DDEs, for which theory and numerical methods are well developed; see, for example, the textbooks \cite{hale2013introduction, smithdde, stepan}.

We are interested here in the effects of delayed negative feedback in systems that are subject to external forcing. The interplay between feedback and periodic forcing can lead to rich and complex dynamics, and is a crucial feature of systems encountered in a wide range of fields, including human balancing \cite{insperger2013}, network dynamics \cite{Semenov2016}, control theory \cite{pyragiene2005}, machine tool vibrations \cite{LelkesKalmarNagy2019}, nonlinear optics \cite{sorrentino2015}, and climate science \cite{GZT, GZT2, SuarezSchopf, TSCJ}. More specifically, the negative feedback term often takes the form of a fast switch between two levels of feedback, which is modeled mathematically by a switching function.

Our motivating example of this type of system is the Ghil--Zaliapin--Thompson (GZT) model \cite{GZT} for the El Ni{\~n}o--Southern Oscillation (ENSO). This climate phenomenon is characterized by irregular switching between phases of a warmer-than-average and cooler-than-average equatorial Pacific ocean, known as El Ni{\~n}o  and La Ni{\~n}a events, respectively, which occur in December to February due to seasonal forcing. The GZT model \cite{GZT} follows the so-called delayed action oscillator (DAO) paradigm \cite{SuarezSchopf, BH}, where wave propagation across the Pacific ocean creates a negative feedback loop; it
describes the deviation from the long-term average of the depth of the thermocline, the sharp interface between warm surface water and the cold deep ocean, and takes the form
\begin{equation}
    \label{eq:GZT}
    h^\prime(t)=-\tanh\bigl(\kappa h(t-\tau)\bigr)+c\cos(2\pi t).
\end{equation}
The first term of~\eqref{eq:GZT} provides switching-type negative feedback by the hyperbolic tangent function $\tanh$ with delay $\tau$; here, the `sensitivity' parameter $\kappa$ is generally set to the accepted standard value $\kappa =11$, which corresponds to fast switching. Time is measured in years, and the second term represents the seasonal periodic forcing of amplitide $c$ (relative to that of the feedback) and period $P = 1$. The GZT model has been studied extensively for the case that the delay $\tau$ is constant. In spite of its simple scalar form, this DDE exhibits rich dynamics and has become a standard conceptual ENSO model for studying resonance phenomena, including the breakdown of invariant tori and the transition to chaos via period-doubling cascades within overlapping resonance tongues \cite{GZT,GZT2,KKGZT,KKChenciner,KKstatedep,KKP_TSCJ}.

Taking delays to be constant is a helpful assumption in the first instance. In many applications, however, relevant response times are not actually constant but depend on the evolving state of the system and may vary considerably. This leads to a mathematical description by state-dependent DDEs. Moreover, in feedback processes, such state dependence is often implicit: the delay selects the past state that drives the feedback, but that past state in turn determines the delay \cite{walther2014topics}. Implicitly defined state-dependent delays arise in applications including echo-based soft landing \cite{walther2007softlanding}, position control \cite{walther2002position}, velocity control \cite{BugerMartin2000, BugerMartin2000_2} and conceptual climate modelling \cite{KKstatedep, Clarke}. Indeed, as highlighted in \cite{chao}, delays in climate processes are generally not constant, and various sources of state dependence exist in the ENSO system. Keane \emph{et al.}~\cite{KKstatedep} recently addressed the question what form such state dependence may take, and whether it influences the dynamics of the GZT model~\eqref{eq:GZT} compared to the constant-delay case. Specifically, they studied how the delay of the feedback depends on the thermocline depth, due to the signal travel time to the surface of the ocean, and this led them to introduce the implicitly defined delay
\begin{equation}
\label{eq:tauT}
\tau=\bar{\tau}+\mu h(t)+\eta h(t-\tau).
\end{equation}
Here, $\bar{\tau}$ is the mean delay due to wave propagation alone, which is still assumed to be constant, and there are two types of linear state dependence: the term $\mu h(t)$ depends on the present thermocline depth, while $\eta h(t-\tau)$ depends on the thermocline depth at time $t-\tau$. Hence, the delay $\tau$ is implicitly defined by~\eqref{eq:tauT} and needs to be solved for, together with $h$. To avoid this additional complexity, the authors of \cite{KKstatedep} considered the GZT model~\eqref{eq:GZT} with the simplified state-dependent delay
\begin{equation}
\label{eq:tauS}
\tau=\bar{\tau}+\mu h(t)+\eta h(t-\bar{\tau}),
\end{equation}
where the delay in $h$ is approximated by the mean delay  $\bar{\tau}$. Their numerical continuation for gradually increasing $\eta$ and $\mu$ showed that instantaneous state dependence alone ($\mu>0$ with $\eta=0$) drastically alters the dynamics, whereas delayed state dependence alone ($\mu=0$ with $\eta>0$) does not appear to have an effect. (A meaningful impact of delayed state dependence was only observed in the simultaneous presence of instantaneous state dependence; see \cite{KKstatedep} for details.) This result is surprising because state-dependence is known to create genuinely new dynamical behavior in other DDEs \cite{Calleja2017ResonancePI, Humphries_2022}. The difference is that, in these earlier works, the state dependence is not delayed itself.

In this paper, we address the question whether there is a structural explanation for the apparent `invisibility' of delayed state dependence that was observed in \cite{KKstatedep} for the GZT model. To this end, we consider the limiting case of $\kappa\to\infty$ when the negative feedback term in~\eqref{eq:GZT} becomes an instantaneous switch, described by the \textit{signum function} $\sign $ as the limit of the switching function $\tanh$. While this may seem quite restrictive, it was established in \cite{Tony2026} that the resulting piecewise-smooth GZT model retains practically all of the dynamics observed in~\eqref{eq:GZT} for the standard value $\kappa = 11$. A major advantage of instantaneous switching is that the piecewise-smooth DDE is considerably more amenable to analysis: because the delayed feedback only takes the values $-1$, $0$ or $+1$, specific solutions and bifurcation loci can be constructed explicitly \cite{Tony2026}. This is the case irrespective of the exact nature of the forcing function, which allows us to establish our main result for the more general class of DDEs
\begin{equation}
\label{eq:psDF}
h^\prime(t) = -\sign \bigl(h(t-\tau)\bigr) + F(t),
\end{equation}
where $F:\mathbb{R}\to\mathbb{R}$ is a measurable forcing function that is bounded in absolute value by some $c > 0$, which represents the `size' of the forcing. We refer to~\eqref{eq:psDF} as the \emph{piecewise-smooth delayed forced} (psDF) system. It generalizes the piecewise-smooth GZT model studied in \cite{Tony2026} as the special case with $F(t)=c\cos(2\pi t)$. Moreover, the psDF system is prototypical and encompasses the square-wave forced model studied in \cite{sgnmodel} (a further simplification of the GZT model), other forcing functions of interest for ENSO, and relay control systems unrelated to ENSO \cite{fridman1993, Shustin1995}.

Our goal is to study the psDF system~\eqref{eq:psDF} with implicit and delayed state dependence. The linear state-dependence expressed by~\eqref{eq:tauT} with $\mu = 0$ represents the simplest possible form. From a modelling perspective, however, it is desirable to allow for other forms that capture nonlinear effects, such as saturation or asymmetry between El Ni{\~n}o and La Ni{\~n}a phases. This is why we formulate our results for implicit delay of the form
\begin{equation}
    \label{eq:tauPsi}
    \tau = \bar{\tau} + \Psi\bigl(h(t-\tau)\bigr),
  \end{equation}
where the response function $\Psi$ encodes how the past state $h(t-\tau)$ modifies the constant delay $\bar{\tau}$. To ensure that~\eqref{eq:tauPsi} appropriately generalizes the linear case recovered by setting $\Psi(u) = \eta u$, we require that
\begin{enumerate}[label=\textbf{(P\arabic*)}]
\item \label{P1} $\Psi$ is Lipschitz continuous, with Lipschitz constant $\eta>0$, which serves again as the measure of the strength of the state-dependence;
\item \label{P2} $\sign (\Psi(u)) = \sign (u)$, so that the sign of the
state (the thermocline anomaly for the GZT model) is preserved.
\end{enumerate}
    
At each fixed time $t$, equation~\eqref{eq:tauPsi} should be viewed as a fixed-point problem for the unknown value of the delay, which is reproduced when evaluating the state $h$ exactly this delay time ago. When the associated fixed-point map is a contraction, then there is a unique solution $\tau_*$. Motivated by this viewpoint and the previous work in \cite{KKstatedep}, we also study nested approximations of the implicit delay. These are obtained by starting from $\bar{\tau}$ and repeatedly substituting the current approximation into the right-hand side of~\eqref{eq:tauPsi}, which generates the sequence of delays
\begin{align}
    \tau_1 &= \bar{\tau}, \nonumber\\
    \tau_2 &= \bar{\tau} + \Psi\bigl(h(t-\bar{\tau})\bigr), \nonumber\\
    \tau_3 &= \bar{\tau} + \Psi\Bigl(h\bigl(t-\bar{\tau}
             -\Psi\bigl(h(t-\bar{\tau})\bigr)\bigr)\Bigr), \nonumber
\end{align}
and so on, via the recurrence
\begin{equation}
    \label{eq:recch}
    \tau_i = \bar{\tau} + \Psi\bigl(h(t-\tau_{i-1})\bigr).
\end{equation}
Clearly, $\tau_1$ recovers the constant-delay case, and $\tau_2$ is the approximation~\eqref{eq:tauS} studied in~\cite{KKstatedep} when $\Psi(u)=\eta u$
and $\mu = 0$.

We are now able to state our main result.

\begin{theorem}[Dynamic invisibility of delayed state-dependence]
  \label{thm:main-invisibility}
  Consider the psDF model~\eqref{eq:psDF} with~\eqref{eq:tauPsi}, where $F$ is measurable and bounded by $c>0$ and the function $\Psi:\mathbb{R}\to\mathbb{R}$ satisfies conditions~\ref{P1} and~\ref{P2} with Lipschitz constant $\eta$.

For an initial Lipschitz continuous function segment $\varphi\in\mathcal{C}([-2\bar{\tau},0],\mathbb{R})$ with Lipschitz constant $L_\varphi$ and satisfying
\[
  \eta|\varphi(s)|<\bar{\tau} \quad \text{for all}\ s\in[-2\bar{\tau},0],
\]
let $g(t)$ be the (unique) solution with initial history $\varphi$ at time $t_0$ of~\eqref{eq:psDF} with constant delay $\tau = \tau_1 = \bar{\tau}$, and assume that it also satisfies 
\[
\eta|g(t)|<\bar{\tau} \quad \text{for all}\  t \geq t_0.
\]

Let $h_i(t)$ be the solution with initial history $\varphi$ at time $t_0$ of~\eqref{eq:psDF} with nested delay $\tau = \tau_i$, defined by the recurrence~\eqref{eq:recch} with $\tau_1 = \bar{\tau}$. 

Provided that 
\begin{equation}
\label{eq:1+c<1}
\eta(1+c)<1 \quad \text{and} \quad \eta L_\varphi<1,
\end{equation}
we have the following.
\begin{itemize}
\item[(i)]
For any $i \geq 2$, we have $h_i(t)=g(t)\ \text{for all} \ t\geq t_0$.
\item[(ii)]
  The psDF~\eqref{eq:psDF} with implicit delay~\eqref{eq:tauPsi} and initial history $\varphi$ at time $t_0$ has the unique solution $h_*(t) =g(t) \ \text{for all} \ t\geq t_0$ with unique delay function $\displaystyle \tau_*(t) = \lim_{i\to\infty} \tau_i(t)$.
\end{itemize}
\end{theorem}

The proof of Theorem~\ref{thm:main-invisibility} is given in Sections~\ref{sec:iteration} and~\ref{sec:proofmaintheorem}. The key ingredient, presented in Section~\ref{sec:iteration}, is that, for a given function segment $x(t)$ on a compact $t$-interval and given fixed delay $\bar{\sigma}$, the associated nested delays $\sigma_j(t)$ are generated by an iterated map; conditions~\eqref{eq:1+c<1} guarantee that this map has the limit $\sigma_*(t)$ as a unique fixed point. As part of this convergence, we show that the zero sets of the delayed function segments $r_j(t) = x(t-\sigma_j(t))$ are the same for any $j \in \mathbb{N} \cup \{*\}$ and, moreover, their signs agree, that is, $\sign (r_1(t)) = \sign (r_j(t))$. The proof in Section~\ref{sec:proofmaintheorem} then consists mainly of the application of this construction to the psDF model~\eqref{eq:psDF} with $x(t)$ being a local segment of $h_1(t) = g(t)$.

\begin{figure}[t!]
   \includegraphics[scale=1]{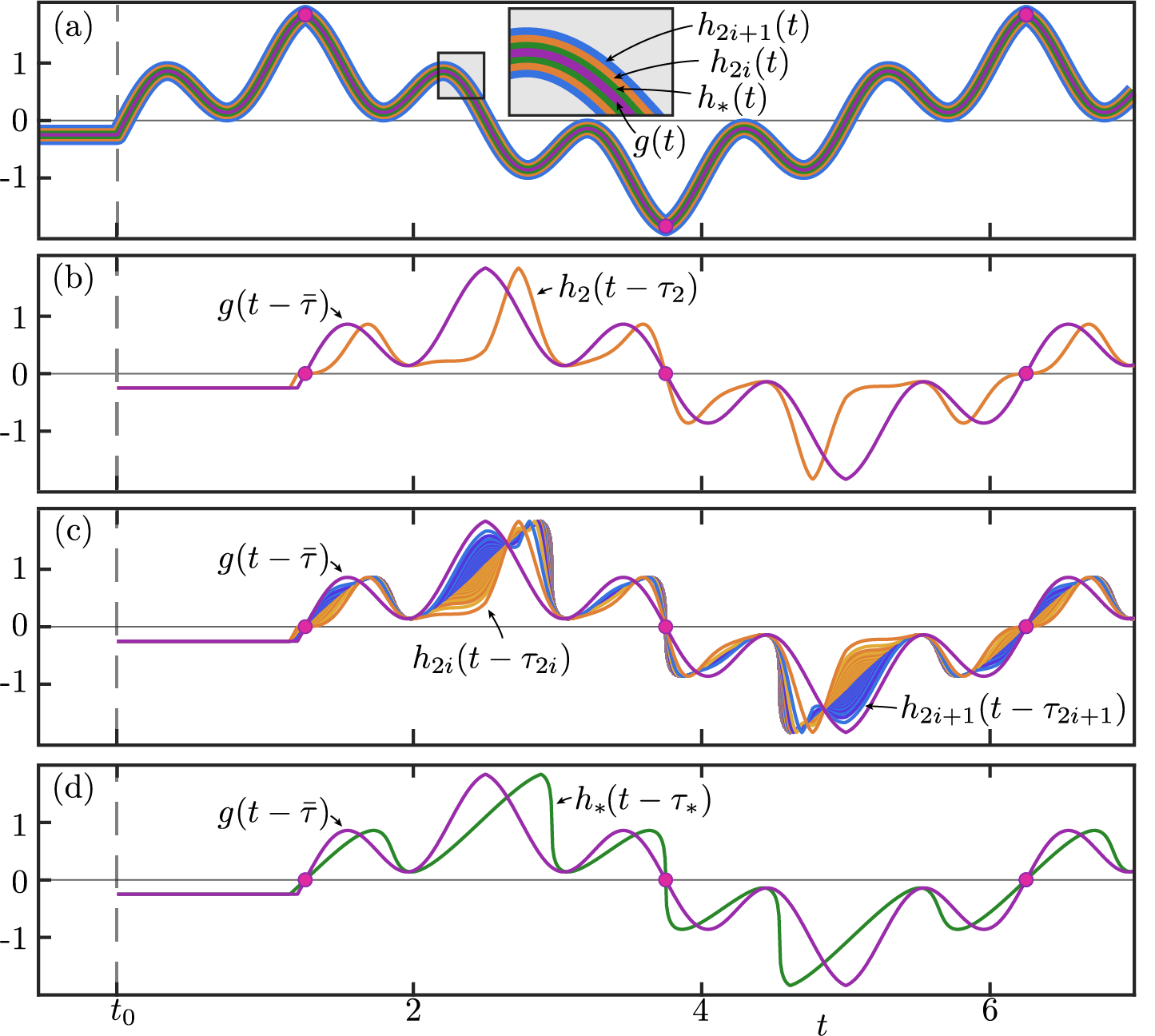}
   \caption{Illustration of Theorem~\ref{thm:main-invisibility} for the piecewise-smooth GZT model, \eqref{eq:psDF} with forcing $F(t)=c\cos(2\pi t)$ and response function $\Psi(u)=\eta u$, for $c=3.7$, $\bar{\tau}=1.25$, $\eta=0.21$, and with initial function $\varphi\equiv-1/4$ at $t_0 \approx -0.054$. Panel~(a) shows the trajectories $g(t)$ with constant delay $\tau = \bar{\tau}$ (purple), $h_i(t)$ with $\tau = \tau_i$ for $ 2 \leq i \leq 150$ (orange for even $i$, blue for odd $i$), and $h_*(t)$ with the implicit delay $\tau_*$ (green). To illustrate that all these trajectories are identical, they are shown with different line thicknesses, as is illustrated further by the enlargement in the inset. The remaining panels compare the delayed term $g(t-\bar{\tau})$ with its state-dependent counterparts, namely: with $h_2(t-\tau_2)$ in panel~(b), with $h_i(t-\tau_i)$ for $ 3 \leq i \leq 150$ in panel~(c), and with $h_*(t-\tau_*)$ in panel~(d). The pink markers indicate jump discontinuities of $g^\prime(t)$ in panel~(a) and zeros of $g(t-\bar{\tau})$ in panels~(b)--(d).}
\label{fig:successTheorem}
\end{figure}

Figure~\ref{fig:successTheorem} illustrates Theorem~\ref{thm:main-invisibility} for the piecewise-smooth GZT (psGZT) model, given by~\eqref{eq:psDF} with forcing $F(t)=c\cos(2\pi t)$ and response function $\Psi(u)=\eta u$; specifically, for constant delay $\bar{\tau}=1.25$, forcing amplitude $c=3.7$, strength $\eta=0.21$ of the state dependence, and initial function $\varphi\equiv-\tfrac{1}{4}$ with $t_0 \approx -0.054$ chosen so that $g(0)=0$. Panel~(a) displays the (unique) trajectory $g(t)$ of the psGZT model for constant delay $\tau = \bar{\tau}$, together with the trajectories $h_i(t)$ for $ 2 \leq i \leq 150$ with nested delays $\tau = \tau_i$ and $h_*(t)$ with implicit delay $\tau = \tau_*$. These trajectories were found with the construction detailed in Appendix~\ref{Appendix1}. For the chosen strength of state dependence we have $\eta(1+c)=0.987$ and $\eta L_\varphi=0$ (since $\varphi$ is constant), so that conditions~\eqref{eq:1+c<1} are satisfied. Indeed, all computed solutions in Figure~\ref{fig:successTheorem}(a) lie exactly on top of each other, that is, are identical, as stated by Theorem~\ref{thm:main-invisibility}. Panels~(b)--(d) show that, due to the different delay functions $\tau_i(t)$, the delayed terms of these solutions are not the same; nevertheless, have the same (nondegenerate) zeros and the same signs in between. Each such zero crossing of the delayed term produces a jump discontinuity in the associated first derivative of the trajectory. It is precisely the sign consistency and the synchronization of zero crossings that explain the agreement observed in panel~(a). Notice the considerable difference already between $g(t-\bar{\tau})$ and $h_2(t-\tau_2)$ in panel~(b). For increasing $i$ the delayed terms $h_i(t-\tau_i)$ converge to $h_*(t-\tau_*)$, as is illustrated in panel~(c); notice here the difference between the delayed terms $h_i(t-\tau_i)$ for even $i$ and those for odd $i$. On the scale of Figure~\ref{fig:successTheorem}, $h_i(t-\tau_i)$ for $i=150$ appears to have converged to the limit $h_*(t-\tau_*)$, which is shown in panel~(d) in comparison with $g(t-\bar{\tau})$.

An interesting and related question is in which order past events are `processed', that is, enter into the evaluation of $h^\prime(t)$ in~\eqref{eq:psDF}. More specifically, the delayed argument $d_i(t)=t-\tau_i(t)$ of the delay function $\tau_i(t)$ records which past time points are used in the evaluation of the relevant delayed term $h_i(d_i(t))$. If $d_i(t)$ is monotonously increasing, then later responses depend on later past events, and the system preserves the classical \emph{temporal order of reactions} (TOR) \cite{walther2009temporalorder}. We will show in Theorem~\ref{thm:TOR} that TOR is satisfied for $\tau_2(t)$ and for the implicit delay $\tau_*(t)$ under conditions~\eqref{eq:1+c<1} of Theorem~\ref{thm:main-invisibility}, while for the intermediate nested delays $\tau_i(t)$ with $3 \leq i < \infty$ TOR requires a stronger condition. When this stronger condition fails while \eqref{eq:1+c<1} still holds, the delayed argument $d_i(t)$ may lose monotonicity --- even though $h_i(t) = g(t)$. In other words, and as we will illustrate, the invisibility of state dependence is a more general property than TOR.

The paper is organized as follows. In Section~\ref{sec:setup} we first consider the psDF model~\eqref{eq:psDF} with constant delay $\tau = \tau_1 = \bar{\tau}$, which serves as the reference throughout; here, we show, for a suitable initial history, that the initial-value problem has a unique global absolutely continuous solution $g(t)$. We also introduce the asymptotic mean condition on the forcing $F(t)$ that ensures that the long-time average of the forcing does not overwhelm the switching slopes $\pm 1$ generated by the delayed feedback. In Section~\ref{sec:iteration} we formulate and prove the main technical lemma that, for a given function segment $x$, the nested delays $\tau_j$ converge to a unique limit $\tau_*(t)$ and, moreover, zeros and signs of the delayed segments $r_j(t)$ all agree. Section~\ref{sec:proofmaintheorem} then presents the proof of Theorem~\ref{thm:main-invisibility}, which boils down effectively to verifying the conditions required for the application of the main lemma to $x(t) = g(t)$.
Section~\ref{sec:TOR} then presents and proves Theorem~\ref{thm:TOR} on the more stringent requirement for TOR; here we also illustrates the failure of TOR for trajectories $h_i(t)$ with nested delays, in spite of them being equal to $g(t)$. In Section~\ref{sec:assymmetric} we then note and formulate as Corollary~\ref{coro:invisibility-asymmetric} that invisibility of state dependence and TOR also apply when the switching function is of a more general asymmetric form, subject to corresponding updates of conditions~\eqref{eq:1+c<1} and~\eqref{eq:condnui}. The final Section~\ref{sec:conclusions} presents an overall conclusion and a brief outlook towards future work, and~\ref{Appendix1} explains how we construct trajectories of~\eqref{eq:psDF} explicitely by tracking the sign changes of the respective delayed term.

\section{Solutions of the psDF model}
\label{sec:setup}

The psDF model~\eqref{eq:psDF} is only piecewise smooth, and we now define what we mean by a solution or trajectory of this state-dependent DDE with either implicit delay~\eqref{eq:tauPsi} or nested delay~\eqref{eq:recch}. State-dependent DDEs are generally studied in a functional framework, in which the delay is also a functional of a solution segment~\cite{stdpgeneral}. In contrast, we treat here any state-dependent delay as a function of time $t$ and work pointwise along the trajectory. This approach is equivalent to the functional framework; namely, whenever the delay is uniquely determined, as is the case under the conditions of Theorem~\ref{thm:main-invisibility}, then this delay, as a function of time, defines a corresponding delay functional of the psDF model~\eqref{eq:psDF} for any solution segment, and vice versa. We employ the pointwise viewpoint because it suits our purposes better: the nested delays $\tau_i(t)$ are defined at each time $t$ by a recurrence, and our analysis in Section~\ref{sec:iteration} compares, at each fixed $t$, the delayed state for state-dependent delay $\tau = \tau_i(t)$ with that for the constant-delay $\tau = \bar{\tau}$.

Our definition of a solution $h_*$ of the psDF model~\eqref{eq:psDF} with implicit delay~\eqref{eq:tauPsi}, which we call the \emph{implicit-delay psDF model}, is inspired by Walther's notion of a solution of an algebraic-delay differential system~\cite{walther2009temporalorder,walther2024solution}, formulated here pointwise in time rather than for delay function segments.

\begin{definition}[Solution of the implicit-delay psDF model]
\label{def:implicit}
Let $t_0\in\mathbb{R}$ and $T\in(t_0,\infty)$. A pair $(h_*,\tau_*)$ is a
\emph{solution} of the implicit-delay psDF model on $[t_0,T)$ with initial
condition $\varphi\in \mathcal{C}([-2\bar{\tau},0],\mathbb{R})$ if
\begin{enumerate}[label=\textbf{(I\arabic*)}]
\item \label{I1} $h_*:[t_0-2\bar{\tau},T)\to\mathbb{R}$ satisfies
$h_*(t)=\varphi(t-t_0)$ for all $t\in[t_0-2\bar{\tau},t_0]$, and $h_*$ is
absolutely continuous on every compact subinterval of $[t_0,T)$;
\item \label{I2}  $\tau_*:[t_0,T)\to(0,2\bar{\tau})$ is continuous;
\item \label{I3}   With $h=h_*$ and $\tau=\tau_*(t)$, the fixed-point equation~\eqref{eq:tauPsi} is satisfied for every $t\in[t_0,T)$, and the DDE~\eqref{eq:psDF} is satisfied for almost every $t\in[t_0,T)$.
\end{enumerate}
\end{definition}

Condition~\ref{I1} is needed because solutions are not continuously differentiable in general: whenever the delayed state $h_*(t-\tau_*(t))$ crosses zero and changes sign, the derivative $h_*^\prime(t)$ jumps; in particular, \eqref{eq:psDF} can be satisfied only almost everywhere. According to~\cite[Section~6.5]{royden}, requiring~\ref{I1} means that $h_*$ is an indefinite integral on $[t_0,T)$; with the requirements in~\ref{I3} it follows that $h_*$ is fully determined by the differential equation. Condition~\ref{I2} is needed because the implicit delay $\tau_*(t)$ may not be unique: without further assumption, \eqref{eq:tauPsi} may have several solutions on some time interval $[t_0,\widehat{T})$, $\widehat{T}<T$, and continuity forces $\tau_*(t)$ to follow a single branch as $t > t_0$ evolves. 

We define in the same spirit the solution $h_i$ of the psDF model~\eqref{eq:psDF} with nested delay~\eqref{eq:recch}, which we call the \emph{nested-delay psDF model} of order $i\in\mathbb{N}$. Here, the nested delays up to $\tau_i$ are determined from $h_i$ by
\begin{equation}\label{eq:tauj-nested}
\tau_1(t)=\bar{\tau},\quad
\tau_j(t)=\bar{\tau}+\Psi\bigl(h_i(t-\tau_{j-1}(t))\bigr),\quad j=2,\dots,i.
\end{equation}

\begin{definition}[Solution of the nested-delay psDF model]
\label{def:nested}
Let $t_0\in\mathbb{R}$, $T\in(t_0,\infty)$ and $i\in\mathbb{N}$. A function
$h_i:[t_0-2\bar{\tau},T)\to\mathbb{R}$ is a \emph{solution} of the
nested-delay psDF model of order $i$ on $[t_0,T)$ with initial condition $\varphi\in \mathcal{C}([-2\bar{\tau},0],\mathbb{R})$ if
\begin{enumerate}[label=\textbf{(N\arabic*)}]
\item \label{N1}  $h_i(t)=\varphi(t-t_0)$ for all $t\in[t_0-2\bar{\tau},t_0]$, and
$h_i$ is absolutely continuous on every compact subinterval of $[t_0,T)$;
\item \label{N2}  for $j=1,\dots,i$ and every $t\in[t_0,T)$, the delay
$\tau_j(t)$ defined by~\eqref{eq:tauj-nested} satisfies
$\tau_j(t)\in(0,2\bar{\tau})$;
\item \label{N3}  With $h=h_i$ and $\tau=\tau_i(t)$, the DDE~\eqref{eq:psDF} is satisfied for almost every $t\in[t_0,T)$.
\end{enumerate}
\end{definition}

In contrast to the implicit-delay case, the nested delays are given explicitly in terms of the solution $h_i$. In particular, the recurrence~\eqref{eq:tauj-nested} can be carried out as long as each delay takes values in $(0,2\bar{\tau})$, so that the next delayed argument $t-\tau_j(t)$ lies in the domain of $h_i$. Note that the DDE is only required to be satisfied (almost everywhere due to jumps of  $h_i^\prime(t)$) for the case that $j=i$, which defines $h_i$.

\subsection{Reference solution $g$ of the constant-delay psDF model}
\label{subsection:constantdelay}

For $i=1$ in Definition~\ref{def:nested}, we recover the constant-delay psDF model, whose solution $g = h_1$ we refer to as the reference solution. Since $0 < \tau_1=\bar{\tau}\in(0,2\bar{\tau})$, condition~\ref{N2} holds automatically in
this case, and only~\ref{N1} and~\ref{N3} remain to be checked to ensure that $g$ is a solution. The DDE~\eqref{eq:psDF} with $\tau =\bar{\tau}$ can be solved by the method of steps \cite{smithdde}: on each interval $I_n = [t_0+(n-1)\bar{\tau}, t_0+n\bar{\tau}]$ with $n\in\mathbb{N}$, the right-hand side of~\eqref{eq:psDF} is a known function of time, and integration yields a unique absolutely continuous solution. 

\begin{proposition}[Reference solutions for constant delay]
  \label{prop:boundness}
For any initial history $\varphi\in\mathcal{C}([-2\bar{\tau},0],\mathbb{R})$, the psDF model~\eqref{eq:psDF} with $\tau = \bar{\tau}>0$, has a unique solution $g$ on $[t_0,\infty)$ in the sense of Definition~\ref{def:nested} with $i=1$. For measurable and bounded $F$, the solution $g(t)$ is bounded on every finite interval. Moreover, if $F$ also satisfies
\begin{equation}
\label{eq:lim-cond}
-1 <\liminf_{t\to\infty}\frac{1}{t}\int_{t_0}^{t}F(u)\,du
\leq \limsup_{t\to\infty}\frac{1}{t}\int_{t_0}^{t}F(u)\,du <1,
\end{equation}
then $g(t)$ cannot be eventually positive or eventually negative.
\end{proposition}

\begin{proof}
The proof of existence directly follows the argument in \cite{Tony2026} for the piecewise-smooth GZT model, which is the special case $F(t) = c \cos(2\pi t)$ of~\eqref{eq:psDF}. By assumption, $g(t)=\varphi(t)$ on $I_0 = [t_0-\bar{\tau}, t_0]$. The right-hand side of~\eqref{eq:psDF} is integrable on $I_1$ and, therefore,
  \[
g(t) = \varphi(0) +
\int_{t_0}^{t} -\sign\!\big(\varphi(u-t_0-\bar{\tau})\big) + F(u)\, du
\]
defines an absolutely continuous solution $g(t)$ on $[t_0-2\bar{\tau},t_0] \cup I_1$; see~\cite[Section~6.5]{royden}. By induction, this yields the absolutely continuous solution $g(t)$ on $I_n$ and, hence, on $[t_0-2\bar{\tau}, t_0 + n \bar{\tau}]$ for any $n$. It follows that $g(t)$ exists for all $t>t_0$ in the sense of Definition~\ref{def:nested} with $i=1$. In particular, $g$ is absolutely continuous on compact subintervals; this also ensures uniqueness, since any two such solutions that are identical to $\varphi(t)$ on $I_0$ must agree for any time $t >0$. 

Finiteness of $g(t)$ on any compact subinterval follows since $|\operatorname{sign}(\cdot)|\le 1$ and
$|F|\le c$. Hence, $|g^\prime(t)|\le 1+c$ almost everywhere and,
for given  $T > t_0$,
\[
  |g(t)|\le|\varphi(0)|+(1+c)(T-t_0)<\infty
  \qquad\text{for all }t\in[t_0,T].
\]

Let now condition~\eqref{eq:lim-cond} hold for $F$, and let us assume that there exists a time $t_1\geq t_0$ such that $g(t)>0$ for all $t\geq t_1$. Then $g(t-\bar{\tau})>0$ for all $t\geq t_1+\bar{\tau}$ and, therefore,
\[
g^\prime(t)=-1+F(t)\quad\text{for almost all} \ t\geq t_1+\bar{\tau}.
\]
Since $g$ is absolutely continuous on $[t_1+\bar{\tau},t]$,
integrating from $t_1+\bar{\tau}$ to $t>t_1+\bar{\tau}$ gives
\[
g(t)=g(t_1+\bar{\tau})-(t-{t_1-\bar{\tau}}) + \int_{t_1+\bar{\tau}}^{t}F(u)\,du
=-t+\int_{t_0}^{t}F(u)\,du+C,
\]
where
\[
C=g(t_1+\bar{\tau})+t_1+\bar{\tau}-\int_{t_0}^{t_1+\bar{\tau}}F(u)\,du
\]
does not depend on $t$. For $t>\max(t_1+\bar{\tau},0)$, dividing by $t$
yields
\[
\frac{g(t)}{t}=-1+\frac1t\int_{t_0}^{t}F(u)\,du+\frac{C}{t}.
\]
Taking the $\limsup$ and using the
strict right inequality in~\eqref{eq:lim-cond} gives
\[
\limsup_{t\to\infty}\frac{g(t)}{t}
=-1+\limsup_{t\to\infty}\frac1t\int_{t_0}^{t}F(u)\,du<0.
\]
Since $t$ is positive, this contradicts $g(t)>0$ for all $t\geq t_1$. 

If instead $g(t)<0$ for all $t\geq t_1$, the anologuous argument gives 
\[
\liminf_{t\to\infty}\frac{g(t)}{t}
=1+\liminf_{t\to\infty}\frac1t\int_{t_0}^{t}F(u)\,du>0,
\]
which is again a contradiction. Hence, $g(t)$ is not eventually positive or negative. 
\end{proof}

Condition~\eqref{eq:lim-cond} expresses the property of the forcing $F(t)$ that the long-time average of the external forcing does not `overwhelm' the switching slopes $\pm1$ of the delayed negative feedback. Several common choices of forcing satisfy this condition:
\begin{itemize}
\item the autonomous case $F\equiv0$;
\item any forcing with $|F(t)|\leq c<1$ almost everywhere  \cite{fridman1993, Shustin1995};
\item any $P$-periodic forcing with mean
\(
\frac{1}{P}\int_0^P F(u)\,du
\)
lying in $(-1,1)$.
\end{itemize}
In particular, the seasonal forcing $F(t)=c\cos(2\pi t)$ of the psGZT model satisfies~\eqref{eq:lim-cond} for every amplitude $c$. A related case is the square-wave forcing $F(t)=c\,\operatorname{sign}\bigl(\sin(2\pi t)\bigr)$ that was introduced in~\cite{sgnmodel}. We remark that the framework we present here also accommodates more detailed representations of forcing relevant to climate modelling. Seasonal forcing may be represented more accurately by finite Fourier series, and one may also consider a forcing term $F(t)=F_{\mathrm{seasons}}(t)+F_{\mathrm{drift}}(t)$, where $F_{\mathrm{seasons}}$ represents seasonal variation and $F_{\mathrm{drift}}$ models a slowly growing monotone drift, for example, due to global warming.

The trajectory $g(t)$ not eventually being of one sign means that $g(t)$ changes sign infinitely many times on $[t_0,\infty)$. In general, accumulation points in the zero set of $g(t)$ are not excluded, and neither are sliding motions (where $g(t)$ becomes temporarily `stuck' at $g=0$). However, when $|F|\leq c<1$, as considered in \cite{fridman1993}, both these phenomena cannot occur \cite{Shustin1995} and the zeros of $g$ are discrete.

\section{Discrete dynamics of nested delays}
\label{sec:iteration}

In this section, we study the nested and implicit delays for a state $x$ in the form of a function segment
\[
x:[-2\bar{\sigma},\Delta]\to\mathbb{R} 
\]
for some fixed $\Delta>0$ and $\bar{\sigma}>0$. We require that
\begin{enumerate}[label=\textbf{(C\arabic*)}]
\item \label{C1}
$x$ is $L$-Lipschitz continuous on $[-2\bar{\sigma},\Delta]$, with
$\eta L<1$;
\item \label{C2}
$\eta\,\|x\|<\bar{\sigma}$,
\end{enumerate}
where $\|x\|$ is the max-norm of $x$ on $[-2\bar{\sigma},\Delta]$ and $\eta$ is the Lipschitz constant of the response function $\Psi$ satisfying conditions~\ref{P1} and~\ref{P2}.

In the proof of Theorem~\ref{thm:main-invisibility}, $x(t)$ will be a local segment of a solution $h_i(t)$ of the psDF model for $i\in\mathbb{N}\cup\{*\}$. However, in the context of this section, we do not require that $x(t)$ is a solution of the psDF model~\eqref{eq:psDF}; to make this clear in the notation, we refer to the fixed function segment as $x$ rather than $h$ and, moreover, denote any associated delays with the letter $\sigma$, rather than $\tau$.

To define nested and implicit delays for $x(t)$, we consider the $t$-family of maps
\begin{equation}
    \label{eq:F}
    \mathcal{F}(\sigma;t):=\bar{\sigma}+\Psi\bigl(x(t-\sigma)\bigr),
    \quad \sigma\in(0,2\bar{\sigma}),\ t\in[0,\Delta].
\end{equation}
For $t\in[0,\Delta]$, the nested delays are then given as iterates $\sigma_j(t) = \mathcal{F}^{j-1}({\bar{\sigma}};t)$ by the recurrence
\begin{equation}
    \label{eq:recurrence}
    \sigma_1(t)=\bar{\sigma},\quad
    \sigma_j(t)=\mathcal{F}\bigl(\sigma_{j-1}(t);t\bigr),\quad j\geq2, 
  \end{equation}
and the implicit delay $\sigma_*(t)$ is a solution of the fixed-point equation
\begin{equation}
    \label{eq:sigma-fixed}
    \sigma_*(t)=\mathcal{F}\bigl(\sigma_*(t);t\bigr),\quad t\in[0,\Delta].
\end{equation}
For $j\in\mathbb{N}\cup\{*\}$, we also define the corresponding delayed states
$r_j:[0,\Delta]\to\mathbb{R}$ of $x(t)$ by
\[
r_j(t):=x\bigl(t-\sigma_j(t)\bigr),\quad t\in[0,\Delta].
\]
Note that $r_1(t)=x(t-\bar\sigma)$ is just a shifted $x(t)$ for $t\in[0,\Delta]$. 

The recurrence~\eqref{eq:recurrence} is a more general and abstract form of the nesting process for the $\tau_i$ that determine a solution $h_i$ of the psDF model~\eqref{eq:psDF}. The key idea is that, for any fixed segment $x$, the associated delays $\sigma_j$ all exist and converge to $\sigma_*$ under suitable conditions, which allows us to conclude that the delayed states $r_i$ all have the same sign. This property is at the heart of the proof of Theorem~\ref{thm:main-invisibility} in Section~\ref{sec:proofmaintheorem}, and we now formulate it as follows.

\begin{lemma}[Implicit delay and sign equality of delayed states]
  \label{lemma:sign}
Suppose that the state $x$ satisfies~\ref{C1}--\ref{C2}. Then the following hold.
\begin{enumerate}
\item[(i)] For every $t\in[0,\Delta]$, the map $\mathcal{F}(\cdot\,;t)$ has a
unique fixed point $\sigma_*(t)$ in $(0,2\bar{\sigma})$, and
$\sigma_j\to\sigma_*$ uniformly on $[0,\Delta]$.
\item[(ii)] The delays $\sigma_j(t)$ for every $j\geq1$ and the delay $\sigma_*(t)$ are Lipschitz continuous on $[0,\Delta]$, with Lipschitz constants
\begin{equation}
  S_j:=\sum_{k=1}^{j-1}(\eta L)^{k}\quad
  \text{and}\quad
  S_*:=\frac{\eta L}{1-\eta L},
  \label{eq:LipschitzSi}
\end{equation}
respectively.
\item[(iii)] For every $t\in[0,\Delta]$ and $j\geq1$,
\[
|r_j(t)-r_*(t)|\leq (\eta L)^{\,j}\,|r_*(t)|\leq (\eta L)^{\,j}\,\|x\|.
\]
\item[(iv)] $r_j\to r_*$ uniformly on $[0,\Delta]$, and for every
$t\in[0,\Delta]$ and $j\geq1$,
\[
\sign\bigl(r_j(t)\bigr)=\sign\bigl(r_*(t)\bigr).
\]
\end{enumerate}
\end{lemma}

\begin{proof}
We consider fixed $t\in[0,\Delta]$, unless stated otherwise, and introduce for notational convenience
\[
M:=\eta\|x\|<\bar{\sigma}\quad \text{and} \quad
X:=[\bar{\sigma}-M,\ \bar{\sigma}+M] \subset (0,2\bar{\sigma}),
\]
whose properties are due to~\ref{C1} and~\ref{C2}.

\smallskip
\noindent(i) The strategy is to show that $\mathcal{F}(\cdot\,;t)$ is a contraction that maps $X$ strictly into itself, so that it has a unique fixed point $\sigma_*(t)$ uniformly in $t \in [0,\Delta]$.

By construction, $x(t-\sigma)$ is well defined for all $\sigma\in(0,2\bar{\sigma})$, and from property~\ref{P1} and $\Psi(0)=0$, which follows from~\ref{P2}, we have
\[
|\Psi(x(t-\sigma))|=|\Psi(x(t-\sigma))-\Psi(0)|
\leq\eta\,|x(t-\sigma)|\leq\eta\|x\|=M,
\]
so that $\mathcal{F}(\sigma;t)=\bar{\sigma}+\Psi(x(t-\sigma))\in X$. In particular, this implies that $\mathcal{F}(\sigma;t)\in X$ for all $\sigma\in X$.

Since $\sigma_1(t)=\bar{\sigma}\in X$ and the $\sigma_j(t)$ are defined by~\eqref{eq:recurrence} as iterates of $\mathcal{F}(\sigma;t)$, we conclude that $\sigma_j(t)\in X$ and, hence, $t-\sigma_j(t) \in (-2\bar{\sigma},\Delta)$. Therefore, all nested delays $\sigma_j(t)$ and corresponding delayed states $r_j(t)$ are well defined.

Let $\sigma,\hat{\sigma}\in(0,2\bar{\sigma})$. Since $\Psi$ is $\eta$-Lipschitz and $x$ is $L$-Lipschitz on $[-2\bar{\sigma},\Delta]$,
\[
\begin{aligned}
|\mathcal{F}(\sigma;t)-\mathcal{F}(\hat{\sigma};t)|
&=|\Psi(x(t-\sigma))-\Psi(x(t-\hat{\sigma}))|\\
&\leq\eta\,|x(t-\sigma)-x(t-\hat{\sigma})|\\
&\leq\eta L\,|(t-\sigma)-(t-\hat{\sigma})|\\
&= \eta L\,|\sigma-\hat{\sigma}|.
\end{aligned}
\]
Since $\eta L<1$ by~\ref{C1}, $\mathcal{F}(\cdot\,;t)$ is a contraction on
$(0,2\bar{\sigma})$ and, hence, also on $X$. Since the interval
$X\subset\mathbb{R}$ is complete, it follows from the contraction mapping theorem that $\mathcal{F}(\cdot\,;t)$ has a unique fixed point $\sigma_*(t)\in X$. Importantly, $\sigma_*(t)$ is also the unique fixed point in all of $(0,2\bar{\sigma})$ because $\mathcal{F}(\cdot\,;t)$ maps $(0,2\bar{\sigma})$ into $X$. We conclude that the implicit delay $\sigma_*(t)$ is uniquely defined for every $t\in[0,\Delta]$, and so is the implicit delayed state $r_*(t)$.

The contraction property of $\mathcal{F}(\cdot\,;t)$ also yields an explicit rate at which the nested delays $\sigma_j(t)$  approach the implicit delay $\sigma_*(t)$. With~\eqref{eq:recurrence} and~\eqref{eq:sigma-fixed} we have
\[
|\sigma_{j+1}(t)-\sigma_*(t)|
=|\mathcal{F}(\sigma_j(t);t)-\mathcal{F}(\sigma_*(t);t)|
\leq \eta L\,|\sigma_j(t)-\sigma_*(t)|.
\]
Induction over $ j\geq1$ gives, with~\eqref{eq:sigma-r}, \ref{P1} and~\ref{C2},
\[
  \begin{aligned}
    |\sigma_j(t)-\sigma_*(t)|
    &\leq (\eta L)^{\,j-1}|\sigma_1(t)-\sigma_*(t)|\\
    &= (\eta L)^{\,j-1}|\bar{\sigma}-(\bar{\sigma}+\Psi(r_*(t))|\\
    &= (\eta L)^{\,j-1}|\Psi(r_*(t))|\\
    &\leq (\eta L)^{\,j-1} M.
\end{aligned}
\]
Since this bound does not depend on $t$, the convergence is uniform in $t$ on $[0,\Delta]$.

\smallskip
\noindent(ii) Let $t,\hat{t}\in[0,\Delta]$ and suppose that for $j \geq 1$ we have
\begin{equation}\label{eq:sigmaj-lip}
|\sigma_j(t)-\sigma_j(\hat{t})|\leq S_j\,|t-\hat{t}|
\end{equation}
with Lipschitz constant $S_j$.
The Lipschitz properties of $\Psi$ and $x$ give
\[
\begin{aligned}
|\sigma_{j+1}(t)-\sigma_{j+1}(\hat{t})|
&=\bigl|\Psi\bigl(x(t-\sigma_j(t))\bigr) -\Psi\bigl(x(\hat{t}-\sigma_j(\hat{t}))\bigr)\bigr|\\
&\leq\eta\,\bigl|x(t-\sigma_j(t))-x(\hat{t}-\sigma_j(\hat{t}))\bigr|\\
&\leq \eta L\,\bigl|(t-\hat{t})-(\sigma_j(t)-\sigma_j(\hat{t}))\bigr|\\
&\leq \eta L\,|t-\hat{t}|+\eta L\,|\sigma_j(t)-\sigma_j(\hat{t})|\\
&\leq \eta L\,(1+S_j)\,|t-\hat{t}|.
\end{aligned}
\]
Hence, we have the recursion $S_{j+1} = \eta L\,(1+S_j)$. Noting that $S_1= 0$ because
\[
  |\sigma_1(t)-\sigma_1(\hat{t})| = |\bar{\sigma}-\bar{\sigma}| = 0,
\]
this gives $S_j$ as in \eqref{eq:LipschitzSi}.

Similarly, from the fixed-point equation~\eqref{eq:sigma-fixed} we obtain 
\[
\begin{aligned}
|\sigma_*(t)-\sigma_*(\hat{t})|
&=|\Psi(x(t-\sigma_*(t)))-\Psi(x(\hat{t}-\sigma_*(\hat{t})))|\\
&\leq\eta\,|x(t-\sigma_*(t))-x(\hat{t}-\sigma_*(\hat{t}))|\\
&\leq\eta L\,\bigl|(t-\hat{t})-(\sigma_*(t)-\sigma_*(\hat{t}))\bigr|\\
&\leq \eta L\,|t-\hat{t}|+\eta L\,|\sigma_*(t)-\sigma_*(\hat{t})|.
\end{aligned}
\]
Hence,
\[
  |\sigma_*(t)-\sigma_*(\hat{t})|\leq \frac{\eta L}{1-\eta L}\,|t-\hat{t}|,
\]
so that $\displaystyle S_*$ is as in \eqref{eq:LipschitzSi}; note that it is the limit of the $S_j$.

\smallskip
\noindent(iii) Note that, in terms of the delayed states, \eqref{eq:recurrence} and~\eqref{eq:sigma-fixed} read
\begin{equation}
    \label{eq:sigma-r}
    \sigma_j(t)=\bar{\sigma}+\Psi\bigl(r_{j-1}(t)\bigr) \ \text{for} \  j\geq2
    \quad\text{and}\quad
    \sigma_*(t)=\bar{\sigma}+\Psi\bigl(r_*(t)\bigr).
  \end{equation}
For $j=1$ we have, with $\sigma_1(t)=\bar{\sigma}$ and $\Psi(0)=0$,
\[
\begin{aligned}
|r_1(t)-r_*(t)|
&=|x(t-\bar{\sigma})-x(t-\sigma_*(t))|\\
&\leq L\,|\bar{\sigma}-\sigma_*(t)|
= L\,|\Psi(r_*(t))|
= L\,|\Psi(0)-\Psi(r_*(t))|\\
&\leq \eta L\,|r_*(t)|.
\end{aligned}
\]
Hence, we conclude for $j\geq2$ by induction that
\begin{eqnarray}
|r_j(t)-r_*(t)|
&=&|x(t-\sigma_j(t))-x(t-\sigma_*(t))| \nonumber \\
&\leq& L\,|\sigma_j(t)-\sigma_*(t)|
=L\,|\Psi(r_{j-1}(t))-\Psi(r_*(t))| \nonumber \\
&\leq& \eta L\,|r_{j-1}(t)-r_*(t)|
\leq (\eta L)^{\,j}|r_*(t)|.
\label{eq:rj-estimate}
\end{eqnarray}
The second inequality in (iii) follows from $|r_*(t)|\leq\|x\|$ for all $t\in[0,\Delta]$.

\smallskip
\noindent(iv) Since the bound $(\eta L)^{\,j}\|x\|$ for $|r_j(t)-r_*(t)|$ from (iii) does not depend on $t$, we have $r_j\to r_*$ uniformly on $[0,\Delta]$ because $\eta L<1$.

The estimate~\eqref{eq:rj-estimate} says that $r_j(t)$ lies in the interval centred at $r_*(t)$ with radius $(\eta L)^{\,j}|r_*(t)| \leq |r_*(t)|$. In particular, for $r_*(t)=0$ this gives $|r_j(t)|\leq 0$, so that $r_j(t)=0$. Suppose now that $r_*(t)>0$. Then
\begin{align*}
r_j(t) &= r_*(t)+\bigl(r_j(t)-r_*(t)\bigr) \geq r_*(t)-|r_j(t)-r_*(t)| \\
&\geq r_*(t)-(\eta L)^{\,j}r_*(t) = (1-(\eta L)^{\,j})\,r_*(t) > 0.
\end{align*}
Likewise, for $r_*(t)<0$ we have 
\[
r_j(t)\ \leq\ r_*(t)+|r_j(t)-r_*(t)|\ \leq (1-(\eta L)^{\,j})\,r_*(t)\ <\ 0.
\]
Hence, as functions of $t \in [0,\Delta]$ the $r_j$ and $r_*$ have the same zeros and the same signs, as stated. 
\end{proof}

\section{Proof of Theorem~\ref{thm:main-invisibility}}
\label{sec:proofmaintheorem}

Throughout the proof, we consider fixed $i$, where either $i \in \mathbb{N}$ with $i \geq 2$, or $i=*$. We consider the unique reference solution $g = h_1$ with the constant delay $\bar{\tau}$, which exists for all $t \geq t_0$ according to Proposition~\ref{prop:boundness}. To apply the framework in Section~\ref{sec:iteration}, we assume for notational convenience that $t_0=0$, which can be achieved by the time shift $t\mapsto t-t_0$. For fixed but arbitrary $T>0$, we define
\[
x(t):=g(t),\quad t\in[-2\bar{\tau}, T].
\]
To apply Lemma~\ref{lemma:sign} to this specific $x$, we verify conditions~\ref{C1}--\ref{C2} with $\Delta = T$.

First of all, $g(t)=\varphi(t)$ is $L_\varphi$-Lipschitz on $[-2\bar{\tau},0]$, by the assumption on $\varphi$. Moreover, $g$ is absolutely continuous according to Proposition~\ref{prop:boundness}, and~\eqref{eq:psDF} gives $|g^\prime(t)|\leq 1+c$ for almost every $t$; hence, $g$ is $(1+c)$-Lipschitz on $[0, T]$. We define $L:=\max(L_\varphi,1+c)$ and consider two times $t_-< 0 < t_+$. The triangle inequality gives
\begin{align*}
|g(t_-)-g(t_+)|
&\leq|g(t_-)-g(0)|+|g(0)-g(t_+)|\\
&\leq L_\varphi(0-t_-)+(1+c)(t_+-0)\leq L(t_+-t_-).
\end{align*}
Hence, $x$ is $L$-Lipschitz on $[-2\bar{\tau},T]$ and the two conditions in~\eqref{eq:1+c<1} translate to $\eta L <1$, which is~\ref{C1}.

Secondly, by assumption, $\eta|\varphi(s)|<\bar{\tau}$ for $s\in[-2\bar{\tau},0]$ and $\eta|g(t)|<\bar{\tau}$ for $t\geq 0$, so that $\eta|x(t)|<\bar{\tau}$ for every $t\in[-2\bar{\tau},T]$. Since $x$ is continuous on the compact interval $[-2\bar{\tau}, T]$, the modulus $|x|$ attains its maximum and, hence, $\eta\|x\|<\bar{\tau}$, which is~\ref{C2}.

Applying Lemma~\ref{lemma:sign} to $x=g$ on $[-2\bar{\tau},T]$ with $\bar{\sigma} = \bar{\tau}$ gives, since the recurrences~\eqref{eq:recch} and~\eqref{eq:recurrence} are identical, the nested delays $\tau_i(t)=\sigma_i(t)$ on $[0,T]$ and the unique limit $\tau_*:=\sigma_*$, all of which lie in $(0,2\bar{\tau})$. Hence,
\[
r_1(t)=g(t-\bar{\tau})
\quad\text{and}\quad
r_j(t)=g\bigl(t-\tau_j(t)\bigr),\quad t\in[0,T],
\]
and we conclude from Lemma~\ref{lemma:sign}(iv) with $j=i$ that 
\[
\sign\bigl(r_i(t)\bigr)=\sign\bigl(r_1(t)\bigr).
\]
It follows with Proposition~\ref{prop:boundness} that 
\[
g^\prime(t)=-\sign\bigl(g(t-\bar{\tau})\bigr)+F(t)
=-\sign\bigl(g(t-\tau_i(t))\bigr)+F(t)
\]
for almost every $t\in[0,T]$. Hence, since $g$ is absolutely continuous by Proposition~\ref{prop:boundness} and all delays take values in $(0,2\bar{\tau})$, a solution $h_i$ of~\eqref{eq:psDF} with delay~$\tau_i$ exists for any $i\geq2$ and $i=*$, and $h_i=g$ on $[0,T]$;
moreover, $\displaystyle \tau_*(t)=\lim_{j\to\infty}\tau_j(t)\in(0,2\bar{\tau})$ is the unique solution of~\eqref{eq:tauPsi} computed with respect to $g$, and $\tau_*(t)$ is continuous by Lemma~\ref{lemma:sign}(ii).

To prove uniqueness of $h_i$, suppose that $\hat{h}_i\neq h_i$ is also a solution on $[0,T]$ of~\eqref{eq:psDF} with initial history $\varphi$ and delay $\hat{\tau}_i$, computed from $\hat h_i$; then
\begin{equation}
  \xi:=\inf\{t\in[0,T]:\hat{h}_i(t)\neq h_i(t)\} < T.
  \label{eq:xi}
\end{equation}
Since $\hat{h}_i$ and $h_i$ are continuous and share the same initial history $\varphi$, we have $\hat{h}_i=h_i=g$ on $[-2\bar{\tau},\xi]$, and every right neighborhood $[\xi,\xi+\delta]$ of $\xi$ contains points where $\hat{h}_i\neq h_i$. To apply Lemma~\ref{lemma:sign} to $x = \hat{h}_i$, we need to ensure condition~\ref{C2}. Since $\hat{h}_i=g$ on $[-2\bar\tau,\xi]$ and the bound $\eta|\hat h_i|<\bar{\tau}$ holds up to $\xi$, by continuity this bound is also satisfied on a short right interval beyond it. We can therefore choose a $\delta>0$ such that
\begin{enumerate}[label={(d\arabic*)}]
\item \label{d1} $\delta<\bar{\tau}$;
\item \label{d2} $\xi+\delta<T$;
\item \label{d3} $\eta\,|\hat h_i(t)|<\bar{\tau}$ for all $t\in[\xi,\xi+\delta]$.
\end{enumerate}
From~\ref{d1} it follows that $t-\bar{\tau}\in[-\bar{\tau},\xi)$ and, hence, $\hat h_i(t-\bar{\tau}) = g(t-\bar{\tau})$ for $t\in[\xi,\xi+\delta]$. 
Condition~\ref{C1} holds for $x=\hat h_i$ exactly as for $g$, since $|\hat h_i^\prime|\le1+c$.

As in the existence part above, we have $\hat{\tau}_i=\sigma_i$, and $\hat{\tau}_*=\sigma_*$ by the uniqueness of the fixed point according to
Lemma~\ref{lemma:sign}(i). Lemma~\ref{lemma:sign}(iv) applied to $x=\hat{h}_i$ with $\Delta=\xi+\delta$ therefore gives, for $t\in[\xi,\xi+\delta]$,
\[
\sign\bigl(\hat h_i(t-\hat\tau_i(t))\bigr)
=\sign\bigl(\hat h_i(t-\bar\tau)\bigr)
=\sign\bigl(g(t-\bar\tau)\bigr),
\]
so that $\hat h_i^\prime=g^\prime$ almost everywhere on $[\xi,\xi+\delta]$. Since $\hat h_i(\xi)=g(\xi)$ and both functions are absolutely continuous, it follows that $\hat{h}_i=g=h_i$ also on $[\xi,\xi+\delta]$, which contradicts the definition of $\xi$.

Since $T$ was arbitrary, the psDF model~\eqref{eq:psDF} with delay $\tau_i$ has the unique solution $h_i= g$ on $[t_0,\infty)$ in the sense of Definitions~\ref{def:implicit} and~\ref{def:nested} for $i\geq 2$ and $i=*$. 
\hfill $\square$

\section{Temporal order in the psDF model}
\label{sec:TOR}

Theorem~\ref{thm:main-invisibility} shows that, under the contraction conditions~\eqref{eq:1+c<1}, the psDF model with any nested delay and with implicit delay psDF generate the same trajectory as when the delay is constant. This invisibility of delayed state dependence does, however, not mean that all these trajectories preserve TOR, that is, that the associated delayed argument functions $d_i(t)=t-\tau_i(t)$, $i\in\mathbb{N}\cup\{*\}$ are strictly increasing. The following result shows in which way TOR is a more restrictive notion than invisibility of delayed state dependence.

\begin{theorem}[Temporal order of reactions in the psDF model]
\label{thm:TOR}
Suppose that the conditions of Theorem~\ref{thm:main-invisibility} hold and, as in Section~\ref {sec:proofmaintheorem}, let
\[
L:=\max(L_\varphi, 1+c)
\]
be the overall Lipschitz constant of $h_i=g$, $i\in\mathbb N\cup\{*\}$, on $[t_0-2\bar\tau,\infty)$.
Then the following hold.
\begin{itemize}
\item[(i)] The nested-delay psDF model~\eqref{eq:psDF} of order $i\ge2$ preserves TOR when $S_i<1$, where $S_i$ is the Lipschitz constant of $\tau_i(t)$ from~\eqref{eq:LipschitzSi}. This is equivalent to
\begin{equation}
  \eta(1+c)<\nu_i\quad\text{and}\quad \eta L_\varphi<\nu_i,
  \label{eq:condnui}
\end{equation}
where $\nu_i$ is the unique positive real solution of the polynomial equation
\begin{equation}
  \sum_{k=1}^{i-1}y^k=1.
  \label{eq:cyclotomic}
\end{equation}
\item[(ii)] The implicit-delay psDF model~\eqref{eq:psDF} preserves TOR without any additional condition.
\end{itemize}
\end{theorem}

\begin{proof}
We assume again, without loss of generality, that $t_0=0$ and consider a fixed but arbitrary time $T>0$. By Theorem~\ref{thm:main-invisibility}, all quantities exist and $h_i = g$. Hence, the delayed argument of $h_i$ is $d_i(t)=t-\sigma_i(t)$ on $[0,T]$ for all $i\in\mathbb N\cup\{*\}$.  Let $0 \leq t_-<t_+ \leq T$ be arbitrary.

\smallskip
\noindent(i) By Lemma~\ref{lemma:sign}(ii) for $x(t) = g(t)$ and  $\Delta=T$, the delay $\tau_i = \sigma_i$ is Lipschitz continuous on $[0,T]$ with Lipschitz constant $S_i$. Therefore, 
\[
d_i(t_+)-d_i(t_-)=(t_+-t_-)-\bigl(\sigma_i(t_+)-\sigma_i(t_-)\bigr)
\geq(1-S_i)(t_+-t_-).
\]
When $S_i<1$, the right-hand side is positive, and so $d_i(t_+)>d_i(t_-)$ as required.

For $i\geq2$ and strictly positive $y$, the map $\displaystyle y\mapsto\sum_{k=1}^{i-1}y^k$ is strictly increasing and equals $1$ at $y=\nu_i$. Since $S_i$ is this map evaluated at $y=\eta L > 0$, we have $S_i<1$ if and only if $\eta L<\nu_i$,  which is equivalent to~\eqref{eq:condnui} since $L=\max(L_\varphi,1+c)$. 

\smallskip
\noindent(ii)
For notational convenience, we introduce $D:= d_*(t_+)-d_*(t_-) = (t_+ - t_-) - (\tau_*(t_+)-\tau_*(t_-))$. From the fixed-point equation~\eqref{eq:sigma-fixed}, noting that~\ref{C1} is satisfied with $x=g$ and $\sigma_* = \tau_*$, we have
\[
  \begin{aligned}
    |\tau_*(t_+)-\tau_*(t_-)|
&= |\Psi\bigl(g(d_*(t_+))\bigr)-\Psi\bigl(g(d_*(t_-))\bigr)| \\
& \leq \eta\,\bigl|g(d_*(t_+))-g(d_*(t_-))\bigr| \\
&  \leq \eta L\,|D|.
\end{aligned}
\]
Resolving the absolute value of the left-hand side and adding $(t_+ - t_-) > 0$ gives
\begin{equation}\label{eq:D-bounds}
(t_+-t_-)-\eta L\,|D|\ \leq\ D\ \leq\ (t_+-t_-)+\eta L\,|D|.
\end{equation}
Suppose now that $D\leq0$, so that $|D|=-D$. Since $1-\eta L>0$ according to~\eqref{eq:1+c<1}, the left inequality in
\eqref{eq:D-bounds} then gives
\[
D\geq \frac{t_+-t_-}{1-\eta L}>0,
\]
which is a contradiction. Hence, $D>0$ and 
\eqref{eq:D-bounds} gives the bounds
\[
\frac{t_+-t_-}{1+\eta L}\;\leq\; d_*(t_+)-d_*(t_-)\;\leq\;
\frac{t_+-t_-}{1-\eta L}.
\]
In particular, $d_*(t_+)>d_*(t_-)$, as required.
\end{proof}

We proceed by stating the properties of the bounding numbers $\nu_i$ in Theorem~\ref{thm:TOR}(i). 

\begin{proposition}
\label{lemma:nu}
The sequence $(\nu_i)_{i\geq2} \subset (\frac{1}{2},1]$ is strictly decreasing, with $\nu_2=1$ and ${\displaystyle \lim_{i\to\infty}} \nu_i = \frac{1}{2}$.
\end{proposition}

\begin{proof}
We remark that equation~\eqref{eq:cyclotomic} and the sequence $(\nu_i)$ arise in other contexts, in particular, in the study of generalized Fibonacci sequences \cite{knacci}. Nevertheless, we prefer to give a direct proof, which is elementary.

Clearly, $\nu_2=1$ and we only consider $i \geq 3$. We define the polynomials
\begin{equation}
  \displaystyle P_i(y) = \sum_{k=1}^{i-1}y^k - 1 = y \frac{1 - y^{i-1}}{1-y} - 1 = \frac{y^i-2y+1}{y-1}
  \label{eq:characeq}
\end{equation}
and conclude that $P_i\left(\frac{1}{2}\right)= -2^{-(i-1)}<0$ and $P_i(1) = i -2 > 0$; moreover, the function $P_i(y)$ is strictly increasing in between, because $\displaystyle P_i^\prime(y) = \sum_{k=1}^{i-1} ky^{k-1}>0$ on $(0,1]$. It follows that there exists a unique $\nu_i \in (\tfrac 12,1]$ with $P_i(\nu_i) = 0$; in particular, $\nu_3 < \nu_2 = 1$. Evaluating $P_{i+1}(y)=P_i(y)+y^i$ at $y = \nu_i$ gives
$P_{i+1}\bigl(\nu_i\bigr)=P_i\bigl(\nu_i\bigr)+\nu_i^{\,i}=\nu_i^{\,i}>0,$
which implies, similarly, that $\nu_{i+1} < \nu_i$. Hence, the sequence $(\nu_i)_{i\geq2}$ is strictly decreasing and bounded below. Therefore, with $ 0 < y < 1$, it follows from~\eqref{eq:characeq} that its limit is the solution of $\displaystyle \lim_{i\to\infty} P_i(y) = \frac{1 - 2y}{y-1} = 0$, giving ${\displaystyle \lim_{i\to\infty}} \nu_i = \frac{1}{2}$.
\end{proof}

\begin{table}[t!]
\centering
\renewcommand{\arraystretch}{1.3}
\begin{tabular}{|c|c|c|c|c|c|c|c|c|c|c|}
\hline
$i$ & $2$ & $3$ & $4$ & $5$ & $6$ & $7$ & $8$ & $9$ & $\cdots$ & $\infty$\\
\hline
$\nu_i$ & $1$ & $0.6180$ & $0.5437$ & $0.5188$ & $0.5087$ & $0.5041$
& $0.5020$ & $0.5010$ & $\cdots$ & $1/2$\\
\hline
\end{tabular}
\vspace{2mm}
\caption{The first elements of the sequence $(\nu_i)_{i \geq 2}$, rounded to four decimal places, and its limit.
}
\label{tab:nu}
\end{table}

Table~\ref{tab:nu} lists the first few $\nu_i$ as computed with Newton's method from the right-hand side of~\eqref{eq:characeq} and rounded to four decimal places. Note the fast convergence of the $\nu_i$ to $\frac{1}{2}$. We mention for context that $\frac{1}{\nu_3}$ is the Fibonacci constant: the largest real solution of the characteristic equation $u^2 - u - 1 =0$ of the Fibonacci sequence, which is the golden ratio. Likewise, $\frac{1}{\nu_i}$ for $i \geq 4$ is the corresponding constant that solves the characteristic equation \cite[Eq.~(3)]{knacci} of the generalized Fibonacci sequence $ f(n) = f(n-1) +f(n-2) + \ldots + f(n-(i-1))$; this can be seen by dividing~\eqref{eq:characeq} by $y^{i-1} > 0$ and setting $u = \frac{1}{y}$, which shows that Proposition~\ref{lemma:nu} is equivalent to \cite[Lemma~3.6]{knacci}.

Because $\nu_2=1$, Theorem~\ref{thm:TOR}(i) shows, in particular, that the nested-delay psDF model with delay $\tau = \bar{\tau} + \Psi\bigl(h(t- \bar{\tau})\bigr)$ preserves TOR under the condition $\eta L < 1$ of Theorem~\ref{thm:main-invisibility}. For $i\geq3$, however, when $\nu_i < \eta L < 1$ TOR is no longer guaranteed, while invisibility of the delayed state dependence still holds. Since $\nu_i>\frac 12$ for all $i$, TOR holds for all $i\geq3$ as well when $\eta L\leq\tfrac12$. On the other hand, the `gap' between $\nu_i$ and $1$ is large for $i\geq3$ and also converges quickly to $\frac{1}{2}$; hence, for any $\frac{1}{2} < \eta L < 1$, the $d_i(t)$ should be expected not to preserve TOR for sufficiently large $i$.

In light of this fact, Theorem~\ref{thm:TOR}(ii) is arguably surprising. The Lipschitz constants $S_i$ of the nested delays $\tau_i$ converge to a Lipschitz constant $S_*=\frac{\eta L}{1-\eta L}$ of the implicit delay; hence, the argument of part~(i) would require $S_*<1$, that is, $\displaystyle \eta L<\tfrac12 = \lim_{i\to\infty}\nu_i$. Nevertheless, the implicit-delay psDF model preserves TOR under the considerably weaker condition $\eta L<1$ of Theorem~\ref{thm:main-invisibility}. The reason is that the implicit delay $\tau_*(t)$ is defined by a fixed-point problem, which, as the proof of Theorem~\ref{thm:TOR}(ii) shows, yields a lower bound on the growth of $d_*(t)$ compared to the $d_i(t)$ of truncated nesting.

\begin{figure}[t!]
\includegraphics[scale=1]{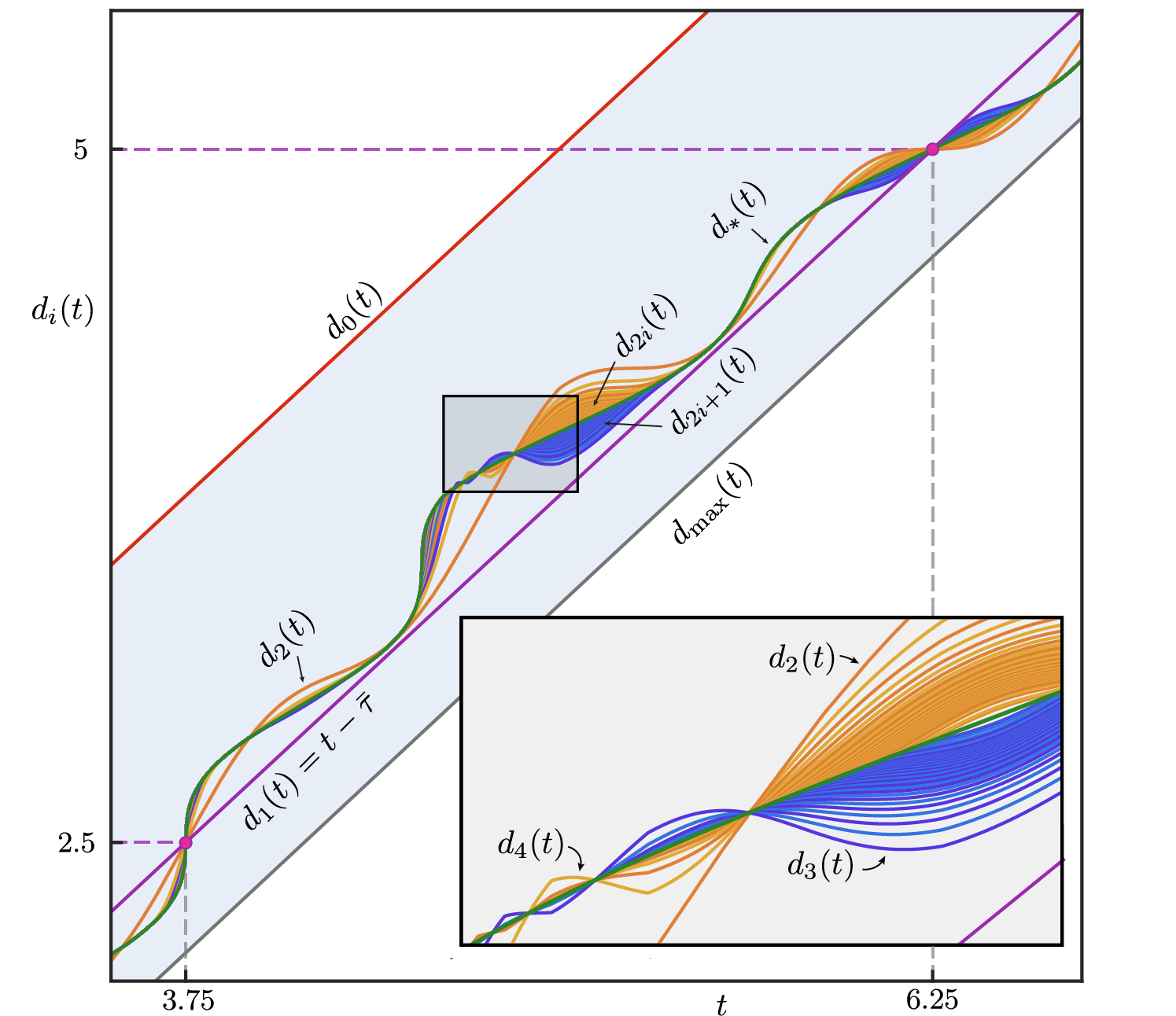}
\vspace{-1mm}
\centering
\caption{
  Illustration of Theorem~\ref{thm:TOR} for the trajectories from Figure~\ref{fig:successTheorem} of the piecewise-smooth GZT model~\eqref{eq:psDF} with forcing $F(t)=c\cos(2\pi t)$ and response function $\Psi(u)=\eta u$. Shown over a smaller range of $t$ are the graphs of the delayed argument functions $d_1(t)=t-\bar{\tau}$ of $g(t)$ (purple), $d_i(t)=t-\tau_i(t)$ of $h_i(t)$ for $2 \leq i \leq 150$ (orange for even $i$, blue for odd $i$), and $d_*(t)=t-\tau_*(t)$ of $h_*(t)$ (green), which all lie between the lines $d_0(t) = t$ of current time (red) and $d_{\text{max}}(t) = t-2\bar{\tau}$ of maximal delay time (grey). The pink markers are at sign changes of the corresponding trajectories $h_i(t) = g(t)$, the dashed horizontal purple lines indicate the times when this happens, and the dashed vertical grey lines indicate when there is the accompanying discontinuity of $g^\prime(t)$. The enlargement in the inset highlights that some the $d_i(t)$ are clearly decreasing on certain time intervals.}
\label{fig:daf1}
\end{figure}

To illustrate Theorem~\ref{thm:TOR}, we consider the trajectories from Figure~\ref{fig:successTheorem} of the piecewise-smooth GZT model~\eqref{eq:psDF} with forcing $F(t)=c\cos(2\pi t)$ and response function $\Psi(u)=\eta u$. Namely, for the parameter values $c=3.7$, $\bar{\tau}=1.25$, $\eta=0.21$ and constant initial function $\varphi$ considered here, we have $\eta L=\eta(1+c)=0.987$ (since $L_\varphi =0$) and, hence,
\[
  \cdots<\nu_4<\nu_3<\eta L =0.987 < \nu_2=1.
\]
Since the condition in Theorem~\ref{thm:TOR}(i) is only sufficient, it is not clear a priori whether TOR actually fails for $i\geq3$. 

Figure~\ref{fig:daf1} shows that the delayed argument functions $d_i$ for $i\geq3$ are indeed decreasing on certain time intervals, while they still converge to the monotonically increasing $d_*$. More specifically, shown are the $d_i(t)=t-\tau_i(t)$ for $1 \leq i \leq 150$, as well as $d_*(t)=t-\tau_*(t)$. Since, by Theorem~\ref{thm:main-invisibility}, $\tau_i(t) \in (0, 2 \bar{\tau})$ for all $i \geq 2$ and $i=*$, the graphs of the $d_i(t)$ all lie in the highlighted diagonal strip $(t - 2 \bar{\tau}, t) =: (d_{\text{max}}(t), d_0(t))$. As Figure~\ref{fig:daf1} and the enlargement in the inset show, $d_2(t)$ is monotonically increasing, but certainly $d_3(t)$ up to $d_9(t)$ can be seen to have $t$-intervals where they are decreasing. Nevertheless, the $d_i(t)$ clearly converge to the monotonically increasing limit $d_*(t)$; namely, from either side depending on whether $i$ is even or odd, as is indicated by alternating colors.

Note that all $d_i(t)$ cross the line $d_1(t) = t - \bar{\tau}$ exactly at the common zeros of the delayed $r_i(t) = h_i(d_i(t)) = g(d_i(t))$; specifically, at $(t, d_i(t)) = (3.75,2.5)$ and at $(t, d_i(t)) = (6.25,5)$ in Figure~\ref{fig:daf1}. Hence, not only do all trajectories $h_i(t) = g(t)$ coincide as is illustrated in Figure~\ref{fig:successTheorem}, but we also see that $d_i(t)=2.5$ and $d_i(t)=5$ all have the unique solutions $t = 3.75$ and $t = 6.25$.  Therefore, even though the $d_i$ may not be monotone, the past events that matter, namely the zeros of the trajectory with a change of sign, are processed at the same times and in the same order for all $i \geq 2$ and $i=*$. Notice further that all $d_i(t)$ are actually monotonously increasing near the crossing points $(t, d_i(t)) = (3.75,2.5)$ and $(t, d_i(t)) = (6.25,5)$. Near their crossing points that do not lie on $d_1(t) = t - \bar{\tau}$, on the other hand, the $d_i(t)$ may be decreasing.

\subsection{The case of an asymmetric switching function}
\label{sec:assymmetric}

To keep the exposition focused, we formulated, proved and illustrated Theorems~\ref{thm:main-invisibility} and~\ref{thm:TOR} for the special case of the symmetric switching function $\sign(\cdot)$. For completeness, we remark that these results can be generalized in a straightforward way to the case that the switching function is asymmetric, that is, the restoring effect in response to positive and negative anomalies differs in magnitude. This is represented by the switching function 
\begin{equation}
A(x)=
\begin{cases}
 a, & x<0,\\
 0, & x=0,\\
 b, & x>0,
\end{cases}
\label{eq:defA}
\end{equation}
for given lower bound $a<0$ and upper bound $b>0$, with $A(x) = \sign(x)$ when $a = -1$ and $b=1$. 

\begin{corollary}[Dynamic invisibility and TOR for asymmetric switching]
\label{coro:invisibility-asymmetric}
The statements of Theorem~\ref{thm:main-invisibility} and Theorem~\ref{thm:TOR} also hold for the psDF model~\eqref{eq:psDF} with switching function $A(\cdot)$ given by~\eqref{eq:defA} instead of $\sign(\cdot)$, subject to 
\begin{equation}
\label{eq:asymmetric_condition}
\eta\bigl(\max(|a|,|b|)+c\bigr)<1
\quad \text{and} \quad
\eta L_\varphi<1
\end{equation}
instead of~\eqref{eq:1+c<1}, and 
\begin{equation}
\label{eq:asymmetric_condition_TOR}
\eta\bigl(\max(|a|,|b|)+c\bigr)<\nu_i
\quad \text{and} \quad
\eta L_\varphi<\nu_i
\end{equation}
instead of~\eqref{eq:condnui}, respectively.
\end{corollary}

All arguments in the proof of Theorem~\ref{thm:main-invisibility} above go through for $A(\cdot)$ with $L=\max(L_\varphi,\max(|a|,|b|)+c)$. First of all, the unique reference solution $g$ with the constant-delay $\bar{\tau}$ also exists for asymmetric switching, because the method of steps in the proof of Proposition~\ref{prop:boundness} also applies since $A$ is bounded and depends only on the already known delayed state from the previous step. Moreover, since $|g^\prime|\leq\max(|a|,|b|)+c$, Lemma~\ref{lemma:sign} applies to $x=g$: the zero sets of the $r_i$ still coincide and their signs agree, and $A(\cdot)$ depends only on the sign of its argument and not on the amplitudes $|a|,|b|>0$. Hence, $h_i = g$ exist for all $i \geq 2$ and $i=*$. The argument for the uniqueness of $h_i$ applies in the same way. Likewise, the proof of Theorem~\ref{thm:TOR} remains unchanged with $L =\max(L_\varphi, \max(|a|,|b|)+c)$.

\section{Outlook and Conclusion}
\label{sec:conclusions}

This work was motivated by the following question: what are the effects of more realistic forms of feedback in ENSO where the delay depends on the state?
A specific motivating example is the interaction between the thermocline and the SST in the Pacific Ocean suggested in the study of the GZT DDE model in \cite{KKstatedep}. In particular, it was found there that delayed state dependence (on its own) does not appear to have an effect for the observed behavior.

To study why this is the case, we considered the prototypical psDF model --- a considerable generalization of the psGZT model --- with implicit delay, which is obtained as the limit of a sequence of nested delays. Crucially, in such piecewise-smooth DDEs, the nonlinearity of the right-hand side is given by a step function, and this property enabled the analytical study presented here. Specifically, we reformulated the nesting construction, so that the nested delays are the iterates of a map, and the implicit delay is its fixed point. This allowed us to prove our main result, Theorem~\ref{thm:main-invisibility}: under conditions on the Lipschitz constants of the initial history and the corresponding solution that guarantee a unique implicit delay via a contraction argument, all nested delays and also the limiting implicit delay generate exactly the same trajectory; namely, that of the psDF model with constant delay. In other words, this type of delayed state dependence, whether nested or implicit, is dynamically invisible. This result holds for quite general types of forcing, for a wider class of response functions and also for an asymmetric switching function. 

We also discussed the temporal order of reactions, or TOR, of the psDF model: under the conditions of Theorem~\ref{thm:main-invisibility}, TOR holds for the second-order and implicit delays, whereas higher-order nested delays require stronger conditions to ensure TOR. Since the respective  trajectories coincide even when these stronger conditions fail, the invisibility of state dependence is a property beyond TOR, as we illustrated with the example of the psGZT model. 

There are several directions for future research. One is the study of the non-contracting regime, where the required condition of Theorem~\ref{thm:main-invisibility} is no longer satisfied. Preliminary results
indicate that delayed state dependence is then no longer invisible and can generate genuinely new dynamics. Without the contraction condition, the implicit delay need not be unique, so that several distinct trajectories may solve the implicit-delay psDF model for the same initial history; it would be interesting to understand how these coexisting solutions are organised. For the nested-delay psDF model, the nested delays then need not converge to any implicit delay as the order of nesting increases. Instead, the nesting map may have periodic orbits and, at the level of the DDE, these can induce richer dynamics at associated levels of nesting.

Interestingly, the numerical study in~\cite{KKstatedep} of the smooth GZT model with $\kappa = 11$, which found that the second-order nested delay has little effect on the dynamics, was carried out in a region of parameter space where the relevant conditions are satisfied. Since Theorem~\ref{thm:main-invisibility} indeed states the invisibility of this type of delay for the piecewise-smooth limit of the GZT model for $\kappa \to \infty$, the work presented here offers an explanation of what was observed in~\cite{KKstatedep}. For constant delay, the piecewise-smooth GZT model was shown to suitably `capture' properties of the smooth GZT model with sufficiently large $\kappa$~\cite{Tony2026}. Therefore, a logical next step will be to investigate how `stable' our main result is when $\kappa$ is finite but large. The proven invisibility of delayed state dependence hinges on the fact that the delayed argument switches between only two distinct slopes, one positive and one negative. Since this property is preserved in good approximation also for large $\kappa$, away from a narrow switching region near $0$, we have the expectation that invisibility is lost in a `gradual way' when $\kappa$ is decreased from large values. The investigation of this process in the smooth GZT model is, hence, a promising direction for future research; this type of study will need to rely on numerical continuation methods, at least, in the first instance.

\section*{Acknowledgments}    
We thank Sam Doak and Kyoung Hyun Lee for insightful discussions. 
This research was supported by Royal Society Te Ap\={a}rangi Marsden Fund grant \#19-UOA-223.



\appendix
\section{Explicit construction of solutions}
\label{Appendix1}

To obtain trajectories, such as those shown in Figure~\ref{fig:successTheorem}, we take advantage of the fact that the psDF model~\eqref{eq:psDF} can be explicitly integrated over intervals on which the sign of the delayed term is constant. We construct solutions in this way instead of using traditional DDE solvers, which generally assume more smoothness than the psDF model offers.

Specifically, we track the sign changes of the delayed term $h_i(t-\tau_i(t))$, which occur at isolated time points for the psGZT model; we refer to them as the switching times and note that each switching time induces a jump discontinuity in $h_i^\prime(t)$. 
Between two consecutive switching times $s_1 < s_2$, we have
\[
h_i(t)=h_i(s_1)-\gamma\,(t-s_1)+\int_{s_1}^t F(u)\,du,
\quad t\in(s_1,s_2),
\]
where $\gamma\in\{-1,1\}$ is the fixed sign of the delayed term on $(s_1,s_2)$. For the psGZT model considered here, $\int_{s_1}^t F(u)\,du=\int_{s_1}^t c\cos(2\pi u)\,du=\tfrac{c}{2\pi}\bigl(\sin(2\pi t)-\sin(2\pi s_1)\bigr)$.

Determining the switching times is straightforward for the constant-delay case: each zero of $h_1=g$ at $t=s$ with a change of sign produces a switching time at $t=s+\bar{\tau}$. For nested delays, the switching times are the times $t$ at which $d_i(t)=t-\tau_i(t) = s$; here, one must also check the genericity condition that $d_i(t)$ actually crosses the level $s$ rather than just touching it. For the implicit delay, we additionally need to solve the fixed-point equation
$\tau_*=\bar{\tau}+\Psi\bigl(h_*(t-\tau_*)\bigr)$ with a Newton-type method. This requires some care, because $t\mapsto h_*(t-\tau_*(t))$ is not differentiable at times $t$ for which $t-\tau_*(t)$ is itself a switching time of $h_*$. We, therefore, find $\tau_*(t)$ by pseudo-arclength continuation in $t$, and check for the switching times separately.

\bibliographystyle{plain}  %
   \bibliography{BSK_implicitSD_references}

@article{SuarezSchopf,
    author =       "Max J. Suarez and Paul S. Schopf",
    title =        "A Delayed Action Oscillator for {ENSO}",
    journal =      "Journal of the Atmospheric Sciences",
    volume =       "45",
    number =       "",
    pages =        "3283–3287",
    year =         "1988",
    DOI =          "10.1175/1520-0469(1988)045<3283:ADAOFE>2.0.CO;2"
}

@article{Humphries_2022,
   title={Nonlinear effects of instantaneous and delayed state dependence in a delayed feedback loop},
   volume={27},
   DOI={10.3934/dcdsb.2022042},
   number={12},
   journal={Discrete and Continuous Dynamical Systems - B},
   publisher={American Institute of Mathematical Sciences (AIMS)},
   author={Humphries, Antony R. and Krauskopf, Bernd and Ruschel, Stefan and Sieber, Jan},
   year={2022},
   pages={7245} }

@article{BugerMartin2000_2,
author = {Matthias Büger and Marcus R.W. Martin},
year = {2000},
pages = {},
title = {The Escaping Disaster: A Problem Related to State-Dependent Delays},
volume = {55},
journal = {Zeitschrift für angewandte Mathematik und Physik ZAMP},
doi = {10.1007/s00033-004-005-6}
}

@article{Calleja2017ResonancePI,
author = {Calleja, Renato and Humphries, Antony R. and Krauskopf, Bernd},
year = {2017},
pages = {},
title = {Resonance Phenomena in a Scalar Delay Differential Equation with Two State-Dependent Delays},
volume = {16},
journal = {SIAM Journal on Applied Dynamical Systems},
doi = {10.1137/16M1087655}
}

@article{walther2024solution,
  title={On the Solution Manifolds for Algebraic-Delay Systems},
  author={Walther, Hans-Otto},
  journal={Ukrainian Mathematical Journal},
  volume={75},
  number={12},
  pages={1801--1821},
  year={2024},
  publisher={Springer}
}

@article{GZT,
author = {Ghil, Michael and Zaliapin, Ilya and Thompson, Skip},
year = {2008},
pages = {},
title = {A delay differential model of {ENSO} variability: Parametric instability and the distribution of extremes},
volume = {15},
journal = {Nonlinear Processes in Geophysics},
doi = {10.5194/npg-15-417-2008}
}

@article{sorrentino2015,
       author = {{Sorrentino}, Taciano and {Quintero-Quiroz}, Carlos Alberto and {Aragoneses}, Andr{\'e}s and {Torrent}, Maria del Carmen and {Masoller}, Cristina},
        title = "{Effects of periodic forcing on the temporally correlated spikes of a semiconductor laser with feedback}",
      journal = {Optics Express},
         year = {2015},
       volume = {23},
       number = {5},
        pages = {5571},
          doi = {10.1364/OE.23.005571},
       adsurl = {https://ui.adsabs.harvard.edu/abs/2015OExpr..23.5571S}
}

@article{pyragiene2005,
author = {Pyragiene, Tatjana and Pyragas, Kestutis},
year = {2005},
pages = {026203},
title = {Delayed feedback control of forced self-sustained oscillations},
volume = {72},
journal = {Physical review. E, Statistical, nonlinear, and soft matter physics},
doi = {10.1103/PhysRevE.72.026203}
}

@article{Semenov2016,
author = {Semenov, Vladimir and Zakharova, Anna and Maistrenko, Y. and Schöll, Eckehard},
year = {2016},
pages = {},
title = {Delayed-feedback chimera states: Forced multiclusters and stochastic resonance},
volume = {115},
journal = {EPL (Europhysics Letters)},
doi = {10.1209/0295-5075/115/10005}
}

@article{insperger2013,
author = {Insperger, Tamas and Milton, John and Stépán, Gábor},
year = {2013},
pages = {20120763},
title = {Acceleration feedback improves balancing against reflex delay},
volume = {10},
journal = {Journal of the Royal Society, Interface / the Royal Society},
doi = {10.1098/rsif.2012.0763}
}

@article{Shustin1995,
author = {Shustin, Eugenii},
year = {1995},
pages = {199-219},
title = {Super-high-frequency oscillations in a discontinuous dynamic system with time delay},
volume = {90},
journal = {Israel J. Math.},
doi = {10.1007/BF02783213}
}

@article{walther2007softlanding,
author = {Walther, Hans-Otto},
year = {2007},
pages = {593-622},
title = {On a Model for Soft Landing with State-Dependent Delay},
volume = {19},
journal = {Journal of Dynamics and Differential Equations},
doi = {10.1007/s10884-006-9064-8}
}

@article{walther2002position,
author = {Walther, Hans-Otto},
year = {2002},
pages = {},
title = {Stable periodic motion of a system with state dependent delay},
volume = {15},
journal = {Differential and Integral Equations},
doi = {10.57262/die/1356060779}
}

@article{walther2009temporalorder,
author = {Walther, Hans-Otto},
year = {2009},
pages = {195-232},
title = {Algebraic-Delay Differential Systems, State-Dependent Delay, and Temporal Order of Reactions},
volume = {21},
journal = {Journal of Dynamics and Differential Equations},
doi = {10.1007/s10884-009-9129-6}
}

@article{fridman1993,
author = {Fridman, Leonid and Fridman, Emilia and Shustin, Eugenii},
year = {1993},
pages = {59-73},
title = {Steady modes in discontinuous control systems with time delay},
volume = {4},
journal = {Pure Mathematics and Applications}
}

@book{hale2013introduction,
  title={Introduction to functional differential equations},
  author={Hale, Jack K. and {Verduyn Lunel}, Sjoerd M.},
  volume={99},
  year={2013},
  publisher={Springer Science \& Business Media}
}

@article{GZT2,
  title={A delay differential model of {ENSO} variability--{P}art 2: {P}hase locking, multiple solutions and dynamics of extrema},
  author={Zaliapin, Ilya and Ghil, Michael},
  journal={Nonlinear Processes in Geophysics},
  volume={17},
  number={2},
  pages={123--135},
  year={2010},
  publisher={Copernicus Publications G{\"o}ttingen, Germany}
}

@article{KKGZT,
    author =       "Andrew Keane and Bernd Krauskopf and Claire Postlethwaite",
    title =        "Delayed Feedback Versus Seasonal Forcing: Resonance Phenomena in an {El {N}i{\~n}o} {S}outhern {O}scillation Model",
    journal =      "SIAM Journal on
Applied Dynamical Systems",
    volume =       "14",
    number =       "",
    pages =        "1229-1257",
    year =         "2015",
    DOI =          "10.1137/140998676"
}

@article{KKChenciner,
    author =       "Andrew Keane and Bernd Krauskopf",
    title =        "Chenciner bubbles and torus break-up in a periodically forced
delay differential equation Model",
    journal =      "Nonlinearity",
    volume =       "31",
    number =       "",
    pages =        "R165",
    year =         "2018",
    DOI =          "10.1088/1361-6544/aab8a2"
}

@article{chao,
    author =       "Andrew Keane and Bernd Krauskopf and Claire M. Postlethwaite",
    title =        "Climate models with delay differential equations",
    journal =      "Chaos",
    volume =       "27",
    number =       "",
    pages =        "114309",
    year =         "2017",
    DOI =          "10.1063/1.5006923"
}

@article{TSCJ,
author = {Eli Tziperman  and Lewi Stone  and Mark A. Cane  and Hans Jarosh },
title = {El {N}i{\~n}o Chaos: Overlapping of Resonances Between the Seasonal Cycle and the Pacific Ocean-Atmosphere Oscillator},
journal = {Science},
volume = {264},
number = {5155},
pages = {72-74},
year = {1994},
doi = {10.1126/science.264.5155.72}}

@article{KKstatedep,
    author =       "Andrew Keane and Bernd Krauskopf and  Henk A. Dijkstra",
    title =        "The effect of state dependence in a delay differential equation model for the {El {N}i{\~n}o {S}outhern {O}scillation}",
    journal =      "Philosophical Transactions: Mathematical, Physical and Engineering Sciences",
    volume =       "377",
    number =       "2153",
    pages   = "20180121",
    year =         "2019",
    DOI =          "10.1098/rsta.2018.0121"
}

@article{KKP_TSCJ,
    author =       "Andrew Keane and Bernd Krauskopf and  Claire M. Postlethwaite",
    title =        "Investigating irregular behavior in a model for the {El {N}i{\~n}o {S}outhern {O}scillation} with positive and negative delayed feedback",
    journal =      "SIAM Journal on
Applied Dynamical Systems",
    volume =       "15",
    number =       "",
    pages =        "1656–1689",
    year =         "2016",
    DOI =          "10.1137/16M1063605"
}

@article{Clarke,
    author =       "Allan J. Clarke and Jianguo Wang and Stephen {van Gorder}",
    title =        "A simple warm-pool displacement {ENSO} model",
    journal =      "Journal of Physical Oceanography",
    volume =       "30",
    number =       "",
    pages =        "1679-1691",
    year =         "2000",
    DOI =          "10.1175/1520-0485(2000)030<1679:ASWPDE>2.0.CO;2"
}

@article{sgnmodel,
    author =       "Pierce Ryan and Andrew Keane and Andreas Amann",
    title =        "Border-collision bifurcations in a driven time-delay system",
    journal = "Chaos",
volume  = "30",
number  = "2",
pages   = "023121",
year    = "2020",
    DOI =          "10.1063/1.5119982"
}

@book{smithdde,
    author    = "Hal Smith",
    title     = "An Introduction to Delay Differential Equations with Sciences Applications to the Life Sciences",
    year      = "2011",
    publisher = "Springer-Verlag",
    address   = ""
}

@article{BH,
  title={Interannual variability in a tropical atmosphere--ocean model: Influence of the basic state, ocean geometry and nonlinearity},
  author={Battisti, David S. and Hirst, Anthony C.},
  journal={Journal of Atmospheric Sciences},
  volume={46},
  number={12},
  pages={1687--1712},
  year={1989}
}

@article{Tony2026,
    author = {Bolduc-St-Aubin, Samuel and Humphries, Antony R.},
    title = {Seasonal forcing dominated dynamics of a piecewise smooth {Ghil-Zaliapin-Thompson} {ENSO} model},
    journal = {Journal of Computational Dynamics},
    volume = {15},
    number = {1},
    pages = {113-162},
    year = {2026},
    doi = {10.3934/jcd.2026011}
}

@incollection{stdpgeneral,
title = {Chapter 5 Functional Differential Equations with State-Dependent Delays: Theory and Applications},
editor = {Ca{\~n}ada, A. and Dr{\'a}bek, P. and Fonda, A.},
booktitle = {Handbook of Differential Equations: Ordinary Differential Equations},
publisher = {North-Holland},
volume = {3},
pages = {435-545},
year = {2006},
issn = {1874-5725},
doi = {10.1016/S1874-5725(06)80009-X},
author = {Ferenc Hartung and Tibor Krisztin and Hans-Otto Walther and Jianhong Wu}
}

@article{walther2014topics,
  title={Topics in delay differential equations},
  author={Walther, Hans-Otto},
  journal={Jahresbericht der Deutschen Mathematiker-Vereinigung},
  volume={116},
  number={2},
  pages={87--114},
  year={2014},
  publisher={Springer}
}

@article{BugerMartin2000,
title = {Stabilizing control for an unbounded
state-dependent delay equation},
journal = {Conference Publications},
volume = {2001},
number = {Special},
pages = {56-65},
year = {2001},
doi = {10.3934/proc.2001.2001.56},
author = {Matthias Büger and Marcus R.W. Martin}
}

@article{LelkesKalmarNagy2019,
author = {Lelkes, J{\'a}nos and Kalm{\'a}r-Nagy, Tam{\'a}s},
title = {Bifurcation analysis of a forced delay equation for machine tool vibrations},
journal = {Nonlinear Dynamics},
volume = {98},
pages = {2961--2974},
year = {2019},
doi = {10.1007/s11071-019-04984-w}
}

@book{stepan,
author = {St{\'e}p{\'a}n, G{\'a}bor},
year = {1989},
title = {Retarded dynamical systems: stability and characteristic functions},
publisher={Longman Scientific \& Technical}
}

@book{royden,
  author    = {Royden, Halsey L. and Fitzpatrick, Patrick M.},
  title     = {Real Analysis},
  edition   = {4th},
  publisher = {Pearson},
  year      = {2010},
  isbn      = {9780131437470},
  url       = {https://www.worldcat.org/title/oclc/456837090}
}

@article{knacci,
author    = {Wolfram, David A.},
title = {Solving Generalized {F}ibonacci Recurrences},
journal = {Fibonacci Quarterly},
volume = {3},
issue = {2},
pages = {129–145},
year = {1998},
doi = {10.1080/00150517.1998.124289489}
}

\end{document}